\documentclass[12pt,reqno]{amsart}
\usepackage{geometry}                
\usepackage{subfig}
\usepackage{amsmath,amssymb}
\usepackage{graphicx}
\usepackage{amssymb,latexsym, amsmath}
\usepackage{amsfonts,amsthm}
\usepackage{amscd}
\usepackage{mathrsfs}
\usepackage{hyperref}
\usepackage{eucal}
\usepackage{multicol}
\usepackage{epstopdf}
\usepackage{amsgen}
\usepackage{xspace}
\usepackage{verbatim}
\usepackage{stmaryrd}
\usepackage{graphicx}
\usepackage{caption}

\newtheorem{theorem}{Theorem}
\newtheorem{proposition}{Proposition}
\newtheorem{lemma}{Lemma}

\newtheorem{remark}{Remark}

\newcommand{\be}{\begin{equation}}
\newcommand{\ee}{\end{equation}}

\newcommand{\ds}{\displaystyle}

\begin{document}
\author[D. Albanez]{Débora A. F. Albanez}
\address{Departamento Acadêmico de Matemática\\ Universidade Tecnológica Federal do Paraná\\ Cornélio Procópio PR, 86300-000, Brazil} \email{deboraalbanez@utfpr.edu.br}

\author[M. Benvenutti]{Maicon José Benvenutti}
\address{Departamento de Matemática\\ Universidade Federal de Santa Catarina\\ Blumenau SC 89065-200, Brazil} 
\email{m.benvenutti@ufsc.br}
\address{Department of Mathematics\\ Trinity Christian College\\ IL 60463, USA} 
\email{mbenvenutti@trnty.edu}

\author[J. Tian]{Jing Tian}
\address{Department of Mathematics \\Towson University\\
Towson, MD 21252, USA} \email{jtian@towson.edu}

\def\bar{\overline}

\title[Recovering the Parameter 
$\alpha$ in the Simplified Bardina Model]
{Recovering the Parameter 
$\alpha$ in the Simplified Bardina Model through Continuous Data Assimilation}

\begin{abstract}
In this study, we develop a continuous data assimilation algorithm to recover the parameter $\alpha$ in the simplified Bardina model. Our method utilizes the observations of finitely many Fourier modes by using a nudging framework that involves recursive parameter updates. We provide a rigorous convergence analysis, showing that the approximate parameter approaches the true value under suitable conditions, while the approximate solution also converges to the true solution.
\end{abstract}
\maketitle
\section{Introduction}


The Navier-Stokes equations (NSE) have been
widely used to describe the motion of viscous incompressible fluid flows. However, solving NSE using the
direct numerical simulation method for turbulent flows is extremely difficult (see \cite{Temam2024}).
Turbulence modeling could provide qualitative and in some cases quantitative measures for a broad spectrum of applications. 
In recent decades, various $\alpha$-regularization models (Navier Stokes-$\alpha$, Leray-$\alpha$, Modified Leray-$\alpha$, Clark-$\alpha$, and the simplified Bardina model) were introduced as efficient
subgrid scale turbulence models (see \cite{bardina1980improved},
\cite{chen1999camassa}, 
\cite{cheskidov2005leray},  \cite{foias2002three}, \cite{IlyinLunasinTiti2006}, and references therein). Derived through an averaging process, these $\alpha$-regularization models not only capture the large scale dynamics of the flow, but also provide reliable closure
models to the averaged equations. Moreover, unlike the other subgrid closure models which usually enhance
dissipation, these $\alpha$-models smooth the nonlinearity of the Navier-Stokes equations. 

The simplified Bardina model is one of the 
$\alpha$-regularization models, which was originally introduced as a closure approximation for the Reynolds equations (\cite{LaytonLewandowski2006}). It has nice analytical, empirical and computational properties, such as global regularity and good 
matching with empirical data collected from turbulent channels and pipes. The explicit steady-state solutions of the simplified Bardina model match the experimental data. Notably, when comparing the number of degrees of freedom in the long-term dynamics of the solutions, the simplified Bardina model has fewer degrees of freedom than the Navier-Stokes-$\alpha$, making it more tractable computationally and analytically \cite{cao2006global}.

This work focuses on the development of a parameter recovery algorithm for the simplified Bardina model. 
\begin{eqnarray}\label{Bardina1}
\left\{
\begin{array}{l}
v_{t}-\nu\Delta v+(u\cdot\nabla)u+\nabla p=f, \\\\
\nabla\cdot v=\nabla\cdot u=0,   \\
\end{array}\right.
\end{eqnarray}
where
\begin{equation}\label{uandv}
v=u-\alpha^2\Delta u.    
\end{equation}
The vector $u=(u_1,u_2,u_3)$ is the spatially filtered velocity field, $p=p(x,t)$ is the modified scalar pressure field, $f=f(x,t)$ is an external force, and $\nu>0$ is the kinematic viscosity of the fluid. In this setting, $f$ and $\nu$ are fixed and known exactly, while the lengthscale parameter $\alpha>0$ and the initial condition $u_0=u(0)$ are assumed to be unknow.

When considering a dynamical system designed to model a physical phenomenon, it is essential that the parameters introduced during the modeling process be carefully chosen to ensure that the model closely reflects the physical reality it is intended to depict, while also contributing to computational efficiency. In applied settings, the determination of these parameters relies on empirical observations. The aim of this work is to develop a recovery algorithm for the parameter $\alpha$, thereby improving the accuracy and reliability of the simplified Bardina model. The parameter $\alpha$ is a key component of the model. It has both a physical and a mathematical origin. Physically, it arises from filtering and averaging, acting as a filter length scale - it specifies the cutoff between large and small scales. Mathematically, it is a regularization parameter which smooths the nonlinear convection term and improves the analytical tractability of the model. It's straightforward to observe that when $\alpha \to 0$, the simplified Bardina model reduces to the NSE.

This work is also inspired by recent studies on parameter recovery using direct observational measurements of the velocity field. In particular, it builds on the Azouani–Olson–Titi (AOT) approach \cite{azouani2014continuous}, especially in the context of the 2D NSE, as in \cite{carlson2020parameter}, where the authors propose an algorithm for dynamically recovering the unknown viscosity $\nu$ of the fluid from data-driven observations of the system.
 Moreover, in \cite{Martinez2022}, convergence analysis for two viscosity update rules is provided: one involving instantaneous evaluation in time and the other relying on averaging in time. Also, a similar algorithm was developed in \cite{Martinez2024}, to determine the external driving force when it is considered to be the unknown parameter, and the convergence of the algorithm is also proved, since sufficiently many modes are observed. Moreover, in \cite{carlsonHudson2022} and \cite{PachevWhiteheadMcQuarrie2022}, multi-parameter recovery in chaotic systems was considered. See also  
\cite{PhysRevFluids.9.054602}, \cite{MartinezMurriWhitehead2025}, 
 \cite{NeweyWhiteheadCarlson2025}, 
 and 
\cite{WangJinHuang2023}, for related results.

Based on argument initially presented in \cite{FoiasProdi} in the context of two dimensional Navier-Stokes equations (NSE), where it was shown that 2D NSE have a finite number of determining modes, i.e., the fact that the long time behavior of the solutions can be determined through of the first $N$ orthogonal Fourier projection $P_N$ onto a finite-dimensional space of bounded linear functionals, we propose here an update algorithm for recovering the lengthscale parameter $\alpha>0$ of the three-dimensional viscous simplified Bardina turbulence model \eqref{Bardina1}. Our goal is to obtain an approximate value for $\alpha$ by recursively running a recovery algorithm, proposed and explicitly presented in formula \eqref{algBeta} in Section 3. In our ideal set-up, we assume that we possess information about the flow field $u$ in terms of a continuous time series $\{P_Nu(t)\}_{t\geq 0}$, for $N>0$, that is, the vector field projection onto the subspace determined by its Fourier modes through wave-number $|k|\leq N$. The key idea lies on considering the feedback control system given by 
\begin{eqnarray}\label{Bardina1assimilated}
\left\{
\begin{array}{ll}
z_{t}-\nu\Delta z+(w\cdot\nabla)w+\nabla p=f -\eta(I-\beta^{2}\Delta) (P_N (w)-P_N (u)), \\\\
\nabla\cdot z=\nabla\cdot w=0,&   \\
\end{array}\right.
\end{eqnarray}
where $I$ is the identity operator, $\eta$ is a positive nudging coefficient, 
\begin{equation}\label{zandv}
z=w-\beta^2\Delta w,   
\end{equation}
and the parameter $\alpha > 0$ from the Bardina system is replaced by a new parameter $\beta=\beta_n > 0$ regularly updated from $\beta_{n-1}$, the value obtained at each time step $n\geq 0$ according to the algorithm in formula \eqref{algBeta}, and applied over the $n$-th time interval $[t_n,t_{n+1}]$. Updates are implemented once certain conditions are satisfied, specified in the main theorem of the work. Details of the update scheme are stated in Section 3.

The approach developed in this work differs from those found in the existing literature. For instance, in \cite{carlson2020parameter}, they develop the parameter recovery algorithm by studying the correlation between the viscosity values $|\nu_2-\nu_1|$ and the difference between the observations $|I_h (u)-I_h (v)|$. While, in \cite{animikh}, the authors extend the definition of the determining map to include viscosity
as an input and recovers it by formulating and solving an optimization problem. Our problem is particularly challenging, as the parameter to be recovered occurs in many terms in the equation.  We overcome this difficulty by establishing several auxiliary lemmas for both the physical system and the data assimilation system, and by designing a recursive update algorithm that plays a central role in the recovery process.

The paper is organized as follows. In Section 2, we establish the classical notation and functional framework in which the auxiliary and main results will be proved. In Section 3, we present the algorithm that reconstructs the parameter $\alpha$. In Section 4, we state the main result of the work, which guarantees the convergence of the proposed algorithm. Also, in Section 5, we provide and prove technical auxiliary lemmas concerning solutions of Bardina system and error estimates involving these solutions and approximation solutions. The main result is proved in Section 6. Finally, Conclusions are presented at the end.



\section{Functional analytic framework and preliminaries} \label{sectionPreliminaries}


We present the mathematical framework relevant to the problems under consideration. The results stated here are standard, with proofs available in \cite{constantin1988navier},  \cite{Foias_Manley_Rosa_Temam_2001}, \cite{RobinsonPierre2003} and \cite{Temam1995}. 

Let $\Omega=[0,L]^3$ denote the three-dimensional torus. We denote by $L^{p} $ the usual three-dimensional Lebesgue vector spaces, and for each $s \in \mathbb{R}$, we define the Hilbert space
\begin{equation}
\dot{H}_{s}= \left\{u(x)=\sum_{K \in \mathbb{Z}^{3}\backslash\{0\}}\hat{u}_{K}e^{2\pi i \frac{K\cdot x}{L}};  \,\, \hat{u}_{K}= \overline{\hat{u}_{-K}}, \,\, \sum_{K \in \mathbb{Z}^{3}\backslash\{0\}} |K|^{2s}|\hat{u}_{K}|^{2}<\infty\right\}, \nonumber
\end{equation}
with the inner product
\begin{equation}
(u,v)_{\dot{H}_{s}}= L^{3}\sum_{K \in \mathbb{Z}^{3}/\{0\}} \left(\frac{2\pi|K|}{L}\right)^{2s} \hat{u}_{K}\cdot\overline{\hat{v}_{K}}, \nonumber
\end{equation}
and the closed subspace 
\begin{equation}
\dot{V}_{s}= \left\{u(x)=\sum_{K \in \mathbb{Z}^{3}\backslash\{0\}}\hat{u}_{K}e^{2\pi i \frac{K\cdot x}{L}};  \,\, \hat{u}_{K}= \overline{\hat{u}_{-K}}, \,\, \hat{u}_{K}\cdot K=0, \,\, \sum_{K \in \mathbb{Z}^{3}\backslash\{0\}} |K|^{2s}|\hat{u}_{K}|^{2}<\infty\right\}, \nonumber
\end{equation}
endowed by the norm
$$\|u\|^2_{\dot{V}_s}=L^3\sum_{K \in \mathbb{Z}^{3}\backslash\{0\}}\left(\dfrac{2\pi|K|}{L}\right)^{2s}|\hat{u}_k|^2. $$
It follows that $\dot{V}_{s_{1}} \subset \dot{V}_{s_{2}}$ if $s_{1} \geq s_{2}$ and $\dot{V}_{-s}$ is the dual of $\dot{V}_{s}$, for all $s\geq 0$.

We denote by $\mathcal{P}:\dot{H}_{s}\rightarrow \dot{V}_{s}$ the classical Helmholtz-Leray orthogonal projection given by 
\begin{equation}
 \mathcal{P}u= \sum_{K \in \mathbb{Z}^{3}\backslash\{0\}} \left(\hat{u}_{K}-\frac{K(\hat{u}_{K} \cdot K)}{|K|^{2}} \right)e^{2\pi i \frac{K\cdot x}{L}}, \nonumber 
\end{equation}
 and $A: \dot{V}_{2s} \rightarrow \dot{V}_{2s-2}$ the operator given by 
\begin{equation}
Au=\sum_{K \in \mathbb{Z}^{3}\backslash\{0\}} \frac{4\pi^{2}|K|^{2}}{L^{2}}\hat{u}_{K}e^{2\pi i \frac{K\cdot x}{L}}. \nonumber
\end{equation}
We have that $Au=-\Delta u = -\mathcal{P} \Delta u=-\Delta \mathcal{P} u$.

We adopt the classical notations  $H= \dot{V}_{0}$, $V= \dot{V}_{1}$, $D(A)= \dot{V}_{2}$, $V^{'}= \dot{V}_{-1}$, $D^{'}= \dot{V}_{-2}$, $\|u\|=\|u\|_{L^{2}}$ and $(u,v)=(u,v)_{L^{2}}$. 
 We have the identities
\begin{align}
\|u\|_{H}= \|u\|, \,\,\,\,\, \|u\|_{V}= \| \nabla u\| \,\,\,\,\, \mbox{and} \,\,\,\,\,\|u\|_{D(A)}= \|A u\|, \nonumber
\end{align}
and the  Poincaré inequalities
 \begin{equation}
 \|u\|^{2}\leq\lambda^{-1}_{1}\|\nabla u\|^{2} \,\,\mbox{  for all } u\in V \mbox{  and  }\,\, \|\nabla u\|^{2}\leq\lambda^{-1}_{1}\|Au\|^{2}\,\, \mbox{ for all  } u\in D(A), \label{Poincare}
\end{equation}
where 
\begin{equation}
\displaystyle \lambda_{1}= \frac{4\pi^{2}}{L^{2}}. \label{lambda}
\end{equation}
We recall several particular cases of the Gagliardo-Nirenberg inequalities:
\begin{equation}
\left\{
\begin{array}{ll}
\| g\|_{L^{3}} \leq c\|g\|^{\frac{1}{2}}\| \nabla g\|^{\frac{1}{2}}, & \forall\, g\in V, \\
\| g\| _{L^{4}} \leq c\|g\|^{\frac{1}{4}}\| \nabla g\|^{\frac{3}{4}}, & \forall\, g\in V, \\
\| g\| _{L^{6}}\leq c\|\nabla  g\|, & \forall\, g\in V,
\end{array}\right.  \label{Gagliardo-Nirenberg}
\end{equation}%
where  $c$  is a dimensionless constant.

Moreover, for each $\alpha>0$, we have
\begin{eqnarray}
(I+\alpha^{2}A)^{-1}u = \sum_{K \in \mathbb{Z}^{3}\backslash\{0\}} \frac{L^{2}}{L^{2}+4\alpha^{2}\pi^{2} |K|^{2}}\hat{u}_{K}e^{2\pi i \frac{K\cdot x}{L}}, \nonumber
\end{eqnarray}
together with the estimates
\begin{align}
\|(I+\alpha^{2}A)^{-1}u\| \leq \|u\| \,\,\,\, \mbox{ and } \,\,\,\, \|(I+\alpha^{2}A)^{-1}u\| \leq \frac{1}{\alpha^{2}}\|u\|_{D^{'}}. \label{inequality-1}
\end{align}

Finally, recalling from classical Fourier analysis, the projection onto low Fourier modes is given by
  \begin{align}
 \displaystyle P_{N}u(x)=P_{N}\left(\sum_{K \in \mathbb{Z}^{3}\backslash\{0\}}\hat{u}_{K}e^{2\pi i \frac{K\cdot x}{L}}\right)= \sum_{0<|K|< N}\hat{u}_{K}e^{2\pi i \frac{K\cdot x}{L}},
 \end{align}
where $\widehat{u}_K$ denotes the Fourier coefficient of $u$ corresponding to the wavenumber $K\in\mathbb{Z}^3$, we have
\begin{align}
\|P_N (\varphi) -\varphi\|^{2}_{L^{2}}\leq \frac{\lambda_{1}^{-1}}{N^{2}}\|\nabla\varphi\|^{2}_{L^{2}}, \label{Ih}
\end{align}
 $$\|P_N (\varphi)\|_{L^{2}}\leq \|\varphi\|_{L^{2}},\ \|\nabla P_N (\varphi)\|_{L^{2}}\leq \|\nabla \varphi\|_{L^{2}},\ P_N A=A P_N.$$

Next, we define the bilinear operator $B:V\times V\rightarrow V^{'}$ by  
\begin{equation}
B(u,v)=\mathcal{P}[(u\cdot\nabla)v], \nonumber
\end{equation}
which also extends to $B:V\times H\rightarrow D^{'}$,  $B:D(A)\times H\rightarrow V^{'}$
and the following property
\begin{equation}\label{zero}
 \langle B(u,v),v\rangle=0,   
\end{equation}
where $\langle \,\cdot\,,\,\cdot\,\rangle$ is denotes the appropriate duality pairing. 

With these definitions, the viscous simplified Bardina system \eqref{Bardina1} can be written as
\begin{eqnarray}\label{Bardina1Leray}
\left\{
\begin{array}{ll}
\ds\frac{dv}{dt}+\nu Av+B(u,u)=f, \\\\
v=u+\alpha^2Au,\\\\
\nabla\cdot v=\nabla\cdot u=0,   \\
\end{array}\right.
\end{eqnarray}
where we assume throughout that the forcing term satisfies $f\in L^{\infty}([0,\infty);H)$.


\section{Overview of the parameter-$\alpha$ recovery algorithm}

With mathematical settings presented in Section \ref{sectionPreliminaries}, we also rewrite the system \eqref{Bardina1assimilated} as 

\begin{eqnarray}\label{Bardina1assimilatedLerayProj}
\left\{
\begin{array}{ll}
\ds\frac{dz}{dt}+\nu Az+B(w,w)=f-\eta(I+\beta^{2}A)\mathcal{P}(P_{N}w-P_{N}u), \\\\
z=w+\beta^2Aw,\\\\
\nabla\cdot z=\nabla\cdot w=0.
\end{array}\right.
\end{eqnarray}

The global well-posedness and stabilization results for the system \eqref{Bardina1assimilatedLerayProj} is established in \cite{albanez2018continuous} for a more general class of observables, with the identity operator $I$ in replacement of $I+\beta^2A$ operator in \eqref{Bardina1assimilatedLerayProj}. Similar results for this operator can be obtained through a straightforward process, which we omit here.

We therefore propose the following algorithm to recovery the parameter $\alpha$, explained in details:

\textbf{Initialization}: Set the initial time $t_1=0$. Choose $w_0$ as an arbitrary initial condition for the system. Let $\beta_1>0$ be the initial guess for the length-scale parameter $\alpha$. We assume $$\beta_1^2\in[\alpha_0^2,\alpha_1^2],$$ where $\alpha_0$ and $\alpha_1$ are positive known lower and upper bounds for the unknown parameter $\alpha$.

\textbf{Step 1}: We have all the information for the first time step: $t_1=0$, $w_1(0)=w_0$.

\textbf{Step 2}: For the second time step
\begin{enumerate}

\item Determine the times $\widehat{t}_1$, $t_2$, and also a projection Fourier natural number $\widetilde{N}_1$ such that \eqref{sdvvv} for $n=1$ holds; 
\item Using these values, compute $\eta_1$ and $N_1$ so that the conditions \eqref{estimates000} to \eqref{cond1} are satisfied with $n=1$, 
\item With the obtained parameters, solve the Bardina system:
\begin{eqnarray}\label{stage1}
\left\{
\begin{array}{ll}
\ds\frac{dz_{1}}{dt}+\nu Az_{1}+B(w_{1},w_{1})=f-\eta_{1}(I+\beta^{2}_{1}A)\mathcal{P}(P_{N_{1}}w_{1}-P_{N_{1}}u), \\\\
z_{1}=w_{1}+\beta_{1}^2Aw_{1},\,\, w_1(0)=w_0, \,\,\, t\in[0,t_2].\\

\end{array}\right.
\end{eqnarray}    
which yields $w_1$, the first iteration value of $w$.
\end{enumerate}

\textbf{Step 3}: 
To update $\beta_2$, apply the formula (with $n=1$):
\begin{align}\label{algBeta} 
\beta_{n+1}^2 &=\beta_{n}^2+\dfrac{1}{\delta_{n}}\int_{\hat{t}_{n}}^{t_{n+1}}\left(\frac{d}{dt}P_{N_{n}}g_{n},P_{N_{n}}(u_{t}+\nu Au)\right)+\beta_{n}^2 \left(\frac{d}{dt}\nabla P_{N_{n}}g_{n},\nabla P_{N_{n}}(u_{t}+\nu Au)\right)\nonumber\\&\nonumber\\&+\nu (\nabla P_{N_{n}}g_{n},\nabla P_{N_{n}}(u_{t}+\nu Au))+\nu\beta_{n}^2 (\nabla AP_{N_{n}}g_{n},\nabla P_{N_{n}}(u_{t}+\nu Au)) 
  \nonumber\\& \nonumber\\
    &+ (P_{N_{n}}B(w_{n},w_{n})- P_{N_{n}}B(w_{n}-P_{N_{n}}g_{n},w_{n}-P_{N_{n}}g_{n}),P_{N_{n}}(u_{t}+\nu Au))\nonumber\\&\nonumber \\&+\eta_{n} (P_{N_{n}}g_{n},P_{N_{n}}(u_{t}+\nu Au))+\beta_{n}^2\eta_{n} (\nabla P_{N_{n}}g_{n},\nabla P_{N_{n}}(u_{t}+\nu Au))\,ds,
    \end{align}
 where $g_n:=w_n-u$ and
\begin{equation} \label{deltaN}
\displaystyle\delta_{n}= \int_{\hat{t}_{n}}^{t_{n+1}}\|\nabla P_{N_{n}}(u_{t}-\nu \Delta u)\|^{2}ds.
\end{equation}

\textbf{Recursive Steps}: 
For each $n\geq2$:
\begin{enumerate}
\item Compute $\beta_{n}$ from $\beta_{n-1}$  by applying formula (\ref{algBeta}).
\item On the interval $[t_n,t_{n+1}]$, solve the system 
 \begin{eqnarray}\label{Bardina1LerayCDA}
\left\{
\begin{array}{ll}
\ds\frac{dz_{n}}{dt}+\nu Az_{n}+B(w_{n},w_{n})=f-\eta_{n}(I+\beta^{2}_{n}A)\mathcal{P}(P_{N_{n}}w_{n}-P_{N_{n}}u), \\\\
z_{n}=w_{n}+\beta_{n}^2Aw_{n},\\
\end{array}\right.
\end{eqnarray}
 with the initial condition $$w_{n}(t_{n})=w_{n-1}(t_{n}).$$
 The values $t_{n+1},\ N_n,\ \eta_n$ are chosen so that conditions  \eqref{sdvvv} and \eqref{estimates000}-\eqref{cond1} are satisfied.
 \end{enumerate}

  \textbf{Final step}: The iteration continues until either the final time $T$ is reached or condition (\ref{deltaN}) degenerates to zero.

\textbf{Remarks:}
 \begin{itemize}

 \item When carrying out the actual simulation, if at each iteration we allow the system to run for a sufficiently long duration such that $\hat{t}_{n}-t_{n}$ is large enough, then $\displaystyle e^{\displaystyle  -\frac{\eta_{n}}{4}(\hat{t}_{n}-t_{n})}$ becomes very small. Consequently, we may expect conditions (\ref{cond3})-(\ref{cond1}) to be satisfied. In this case, $\eta_{n}$ and $N_{n}$ should be adjusted so that (\ref{estimates000})-(\ref{cond2}) hold.
     
\item The left-hand sides of (\ref{estimates0001})-(\ref{cond2}) depend on the quantities $M_{1}$, $M_{2}$, $M_{3}$ and $M_{4}$ given in (\ref{M1})-(\ref{M4}). These quantities include negative exponential time terms that can be neglected in a simulation if it runs for a sufficiently long duration.


\item 
In each step, to successfully implement the update of $\beta$ as described in (\ref{algBeta}), we must find a time $t_{n+1}$ and a projection parameter $N_{n}$ such that the expression in (\ref{deltaN}) is positive; see condition (\ref{sdvvv}). However, there may be cases where this condition is not satisfied; in such instances, the algorithm must be stopped, and no further progress in the approximation process is possible. In this situation, for $t \ge t_{n}$ the physical solution $u(t)$ satisfies $$u_{t}-\nu\Delta u=0,$$ and the simplified Bardina equation (\ref{Bardina1}) formally reduces to \begin{eqnarray}\label{Bardina1ghhgh} \left\{ \begin{array}{l} (u\cdot\nabla)u+\nabla p=f, \\\\ \nabla\cdot u=0. \\ \end{array}\right. \end{eqnarray} It is noteworthy that in (\ref{Bardina1ghhgh}), the unknown lenght-scale $\alpha>0$ no longer influences the dynamics, indicating that there is insufficient information in $u(t)$ to support any further approximation of $\alpha$.



\item To provide some intuition behind formula (\ref{algBeta}): it is obtained by taking the inner product of the difference between systems (\ref{Bardina1}) and (\ref{Bardina1assimilated})   with $P_{N_{n}}(u_{t}+\nu Au)$ (see equations (\ref{eqdifference2})-(\ref{update2}) for details). The formula comes from isolating the term $\alpha - \beta$.
A crucial approximation concerns the nonlinear term $B(u,u)$, specifically, we use $$B(P_{N_{n}}u+w_{n}-P_{N_{n}}w_{n},P_{N_{n}}u+w_{n}-P_{N_{n}}w_{n})=B(w_{n}-P_{N_{n}}g_{n},P_{N_{n}}u-P_{N_{n}}g_{n}),$$
since $u= P_{N_{n}}u+(I-P_{N_{n}})u$ and we expect the high Fourier modes approximation
$(I-P_{N_{n}})u\approx (I-P_{N_{n}})w$, after sufficient runtime. 

Observe that the difference of the nonlinear terms in (\ref{algBeta}) can be expressed as:
\begin{align}\label{eqeq22}
    P_{N_{n}}B(w_{n},&w_{n})-P_{N_{n}}B(w_{n}-P_{N_{n}}g_{n},w_{n}-P_{N_{n}}g_{n}) \nonumber\\
   &=-P_{N_{n}}B(P_{N_{n}}g_{n},P_{N_{n}}g_{n})
  +P_{N_{n}}B(w_{n},P_{N_{n}}g_{n})
    +P_{N_{n}}B(P_{N_{n}}g_{n},w_{n}),
 \end{align}
 which is useful in obtaining the correct estimates. Indeed, it is essential to obtain $g_{n}$ in each step and to get the suitable estimates (see (\ref{abs0101})-(\ref{ult1001.1}) for details).

\end{itemize}

\section{Main Results}

Although $u(t)$ and $\alpha$ are considered to be unknown, we assume that, in addition to the measurements given by $P_{N_{n}}u(t)$ over time, there exist known positive constants $\alpha_{0}$, $\alpha_{1}$, $M_{A}$, $M_{B}$ and $M_{C}$ such that 
\begin{eqnarray}\label{estimates00}
\left\{
\begin{array}{l}
\displaystyle\alpha_{0}\leq \alpha\leq \alpha_{1},\\\\
\displaystyle \|u(0)\|\leq M_{A},\,\, \|\nabla u(0)\|\leq  M_{B},\,\, \|\Delta u(0)\| \leq M_{C}.  

\end{array}\right.
\end{eqnarray}

Define the following functions:
\begin{align}
    M^{2}_1(t)&:= e^{-\nu\lambda_{1}t}\left(M^{2}_{A}+\alpha^{2}_{1}M^{2}_{B}\right)+\frac{1}{\lambda_1^2\nu^2}\cdot \sup\limits_{s\geq 0}\|f(s)\|^2, \label{M1}\\\nonumber\\
M^{2}_2(t)&:=e^{-\nu\lambda_{1}t}\left(M^{2}_{B}+\alpha^{2}_{1}M^{2}_{C}\right)+\frac{2c^{4}}{\alpha^{5}_{0}\nu^{2}\lambda_{1}}e^{-\nu\lambda_{1}t}\left(M^{2}_{A}+\alpha^{2}_{1}M^{2}_{B}\right)^{2}\nonumber\\&+\frac{1}{\nu^{2}\lambda_{1}}\sup_{s\geq 0}\|f(s)\|^2+\frac{2c^{4}}{\alpha_{0}^{5}\nu^{6}\lambda_1^5}\sup_{s\geq 0}\|f(s)\|^4, 
\label{M2}\\\nonumber\\
M_3(t)&:= \frac{\nu}{\alpha_{0}} M_{2}(t)+\frac{c^{2}}{\alpha^{4}_{0}\lambda_{1}^{\frac{3}{4}}}M^{2}_{1}(t)+\sup_{s\geq 0}\|f(s)\|,\label{M3}\\\nonumber\\
M^{2}_{(4,n)}(t,\eta_{n})&:=   \displaystyle e^{\displaystyle-\frac{\eta_{n}}{2}(t-t_{n})}\left(\|w_{n-1}(t_{n})\|^2 + \beta_{n}^2\|\nabla w_{n-1}(t_{n})\|^2\right) \nonumber\\&+\frac{4}{\eta^{2}_{n}}\sup_{r\geq 0}\|f(r)\|^2+ 2\left(2+\frac{\beta_{n}^{2}}{\alpha_{0}^{2}}\right)M^{2}_{1}(t_{n}),\label{M4}\\\nonumber
\end{align}
where $\lambda_{1}$ is given in (\ref{lambda}) and $c$ in (\ref{Gagliardo-Nirenberg}).

We can now state the main result. The theorem provides explicit criteria for selecting $\hat{t}_{n}$, $t_{n+1}$, $\eta_{n}$ and $N_{n}$ at each iteration, ensuring that $\beta_{n}\xrightarrow{n \to \infty} \alpha $  and  $w_{n}(t_{n})\xrightarrow{n \to \infty} u(t_{n}) $ exponentially in $\mathbb{R}$ and $\dot{V}_{1}$, respectively.  

\begin{theorem} \label{theo121}
Consider $u$ solution of \eqref{Bardina1Leray} with initial condition $u(0) \in \dot{V}_{2}$ and $w_{0} \in \dot{V}_{2}$.
Let $0<\varepsilon<\alpha^{2}_{0}$ and $\beta_{1} \in [\alpha_{0},\alpha_{1}]$,
with $\alpha_0,\alpha_1$ satisfying \eqref{estimates00}.  

 Suppose that for each $n \in \mathbb{N}$, there exist $\tilde{N}_{n} \in \mathbb{N}$,  $\hat{t}_{n}\geq 0$, and $t_{n+1}>0$  such that 
\begin{align}\label{sdvvv}
\displaystyle \int_{\hat{t}_{n}}^{t_{n+1}}\|\nabla P_{\tilde{N}_{n}}(u_{t}-\nu \Delta u)\|^{2}\,ds>0,
\end{align}
with $t_{n+1} >\hat{t}_{n}\geq t_{n}$, where $t_{n}$ is the final time from the previous iteration and $t_{1}=0$ for $n=1$.
 Under these assumptions, define 
\begin{align}
\displaystyle \zeta_{n}:=\frac{\displaystyle \int_{\hat{t}_{n}}^{t_{n+1}}\|\nabla P_{\tilde{N}_{n}}(u_{t}(s)-\nu \Delta u(s))\|^{2}ds}{t_{n+1}-\hat{t}_{n}}.\label{zeta}
\end{align}

For each $n\geq 2 $, let $\beta_{n}$ be obtained from the previous iteration via the update algorithm (\ref{algBeta}), where the initial iteration is performed using the given $\beta_{1}$. Moreover, choose $\eta_{n}$ and $N_{n}$ sufficiently large so that conditions \eqref{estimates000}-\eqref{cond1} are satisfied.

\begin{eqnarray}\label{estimates000}
\displaystyle 
 \tilde{N}_{n}\leq N_{n};
\end{eqnarray}

\begin{eqnarray}\label{estimates0002}
\displaystyle 
\frac{\eta_{n}}{N_n^2}\leq \frac{\nu\lambda_1 }{2};
\end{eqnarray}

\begin{eqnarray}\label{estimates0001}
\displaystyle \ds\frac{27c^4M_1^4(t_{n})}{8\nu^{3}\alpha_{0}^4}\leq \eta_{n}; 
\end{eqnarray}

\begin{align}
\max\left\{1,\frac{\sqrt{\epsilon}+\alpha_{1}}{\beta_{n}}\right\}\left(\ds\frac{1}{\nu^{\frac{1}{2}}}M_{3}(t_{n})+\ds\frac{\nu^{\frac{1}{2}}}{\alpha_{0}}M_{2}(t_{n})\right)\frac{16}{\beta_{n}}\leq \eta_{n}^{\frac{1}{2}}\nu\lambda^{\frac{3}{4}}_{1};\label{cond4}
\end{align}

\begin{align}
 \displaystyle 
 \displaystyle \frac{8c^{2}\lambda^{-\frac{1}{4}}_{1}}{\eta^{\frac{1}{2}}_{n}\beta^2_{n}\zeta_{n}^{\frac{1}{2}}N_{n}^{\frac{1}{2}}}\left[ \frac{M_{1}(\hat{t}_{n})}{\alpha_{0}}+\frac{M_{(4,n)}(\hat{t}_{n}, \eta_{n})}{\beta_{n}}\right]\left(\frac{\nu^{\frac{1}{2}}M_{2}(t_{n})}{\alpha_{0}}+\frac{M_{3}(t_{n})}{\nu^{\frac{1}{2}}}\right)\leq \frac{\varepsilon}{4|\alpha^{2}_{1}-\alpha^{2}_{0}|\chi_{1}(n)+4\varepsilon};\label{cond2}
\end{align}

\begin{align}
\max\left\{1,\frac{\sqrt{\epsilon}+\alpha_{1}}{\beta_{n}}\right\}e^{\displaystyle  -\frac{\eta_{n}}{4}(t_{n+1}-t_{n})}\leq \frac{1}{8}
 ;\label{cond3}
\end{align}

\begin{align}\label{cond0}
 \displaystyle \frac{4c^{2}}{\zeta_{n}^{\frac{1}{2}}N_{n}^{\frac{1}{2}}\beta_{n}}\left[ \frac{M_{1}(\hat{t}_{n})}{\alpha_{0}}+\frac{M_{(4,n)}(\hat{t}_{n}, \eta_{n})}{\beta_{n}}\right]
  e^{\displaystyle  -\frac{\eta_{n}}{4}(\hat{t}_{n}-t_{n})}\leq\frac{1}{4\nu\lambda^{\frac{1}{2}}_{1}};
\end{align}
\begin{align}\label{cond1}
 \displaystyle \frac{8c^{2}\lambda^{-\frac{1}{4}}_{1}}{\beta_n\zeta_{n}^{\frac{1}{2}}N_{n}^{\frac{1}{2}}}\left[ \frac{M_{1}(\hat{t}_{n})}{\alpha_{0}}+\frac{M_{(4,n)}(\hat{t}_{n}, \eta_{n})}{\beta_{n}}\right]
  e^{\displaystyle  -\frac{\eta_{n}}{4}(\hat{t}_{n}-t_{n})}\left[\left(1+\frac{\beta_{n}}{\alpha_{0}}\right)M_{1}(t_{n}) +M_{4,n}(t_{n},\eta_{n})\right]\leq\frac{\varepsilon}{2},
\end{align}
where $\chi_{1}(1)=1$ and $\chi_{1}(n)=0$ when $n \geq 2$.

Finally, let $w_{n}$ denote the solution of system (\ref{Bardina1LerayCDA}) on the interval $[t_{n},t_{n+1}]$ with initial condition $w_n(t_{n})=w_{n-1}(t_{n})$ and $w_{1}(0)=w_{0}$ when $n=1$.

Then, for all $n\geq 0$, 

\begin{align}
|\beta^{2}_{n+1}-\alpha^{2} |&\leq \frac{\nu^{-1}\lambda_{1}^{-\frac{3}{4}}(\|g_{1}(0)\|  +\beta_{1}\|\nabla g_{1}(0)\|)+|\beta_{1}^2-\alpha^2|}{2^{n}}, \label{rrr2}
\end{align}
\begin{align}
\|g_{n+1}(t_{n+1})\| + \beta_{n+1}\|\nabla g_{n+1}(t_{n+1})\| &\leq \frac{\|g_{1}(0)\|  +\beta_{1}\|\nabla g_{1}(0)\|+\nu\lambda_{1}^{\frac{3}{4}}|\beta_{1}^2-\alpha^2|}{2^{n}}, \label{rrr3}
\end{align}
where $g_{n}(t):=w_{n}(t)-u(t)$.

\end{theorem}
 \begin{remark} 
In the proof of the above theorem, we first establish the auxiliary estimate
     \begin{align}
|\beta^{2}_{n+1}-\alpha^{2}| &\leq \frac{\varepsilon}{2}+ \frac{\varepsilon}{4^{n}},\,\, \forall\, n \geq 1. \label{rrr1}
\end{align}   
From (\ref{rrr1}) and the fact that $0<\varepsilon<\alpha^{2}_{0}\leq \alpha^{2}$, we conclude the updated $\beta^{2}_{n+1}$ is positive. Indeed, we have
$$\beta^{2}_{n+1}\geq -\frac{\varepsilon}{2}- \frac{\varepsilon}{4^{n}}+\alpha^{2}>  -\frac{\varepsilon}{2}- \frac{\varepsilon}{4^{n}}+\varepsilon>0.$$
 Note that conclusions (\ref{rrr2}) and (\ref{rrr3}) show that $\beta_{n}\xrightarrow{n \to \infty} \alpha $  and  $w_{n}(t_{n})\xrightarrow{n \to \infty} u(t_{n}) $ exponentially in $\mathbb{R}$ and $\dot{V}_{1}$, respectively.  
 \end{remark}
\begin{remark}
Let us emphasize some key elements in the proof of the aforementioned theorem:
\begin{enumerate}
\item After taking the inner product of the difference between systems (\ref{Bardina1}) and (\ref{Bardina1assimilated})   with $P_{N_{n}}(u_{t}+\nu Au)$ (see (\ref{eqdifference2})-(\ref{update2})), the term $\alpha^{2}-\beta_{n+1}^{2}$ is isolated, leading to equality (\ref{ult1001.0aaa});
\item
We estimate (\ref{ult1001.0aaa}) in a straightforward way to obtain the result stated in Proposition \ref{teoBardina25}; see estimate (\ref{dfhbb});
\item We estimate the right-hand side of (\ref{dfhbb}), which involves $\|\nabla w_{n}\|$, $\|\nabla u\|$, and $\|\nabla g_{n}\|$, in terms of a combination of $|\beta^{2}_{n}-\alpha^{2}|$ and $\|g_{n}(t_{n})\|  +\beta_{n}\|\nabla g_{n}(t_{n})\|$ (see estimate (\ref{last1})). In this step, we use Lemma \ref{lem1} and Propositions \ref{teoBardina1} and \ref{teoBardina2}. 
\item We also estimate $\|g_{n+1}(t_{n+1})\|  +\beta_{n+1}\|\nabla g_{n+1}(t_{n})\|$  in terms of an expression involving a combination of $|\beta^{2}_{n}-\alpha^{2}|$ and $\|g_{n}(t_{n})\|  +\beta_{n}\|\nabla g_{n}(t_{n})\|$ - see estimate (\ref{last2}). For this step, we use estimate (\ref{rrr1}) and Proposition \ref{teoBardina2}.
\item The two inequalities (\ref{last1})-(\ref{last2}) form the system (\ref{ullltttmmm}). Then, through a straightforward induction argument, we obtain
(\ref{rrr2}) and (\ref{rrr3}). 
\end{enumerate}
\end{remark}


\section{Auxiliary Estimates}



\begin{lemma}\label{lem1}
 Let $u$ be a solution of \eqref{Bardina1Leray} and $M_{1}(t)$ be as defined in (\ref{M1}), we have
\begin{equation}
\|u(t)\|^2+\alpha^2\|\nabla u(t)\|^2\leq M^{2}_1(t),\,\, \forall\, t \geq 0. \label{ine01}
\end{equation}
\end{lemma} 

The above result is an immediate consequence of Lemma 2 of \cite{ALBANEZ2025109073}.




\begin{lemma}\label{lem12}
For a solution $u$ of \eqref{Bardina1Leray}, with
$M_{2}(t)$ defined in (\ref{M2}),
we have 
\begin{equation}
\|\nabla u(t)\|^2+\alpha^2\|A u(t)\|^2\leq M^{2}_2(t),\,\, \forall\, t \geq 0. \label{ine011}
\end{equation}

\end{lemma} 


\begin{proof}
Multiplying system \eqref{Bardina1Leray} by the solution $Au(t)$, integrating over the domain $\Omega=[0,L]^{3}$, integrating by parts and using (\ref{Poincare}), (\ref{Gagliardo-Nirenberg}) and Young's inequality, we obtain
\begin{align}
    \frac{1}{2}\frac{d}{dt}(\|\nabla u(t)\|^2+&\alpha^2\|A u(t)\|^2)+\nu(\|A u(t)\|^2+\alpha^2\|\nabla Au(t)\|^2)\nonumber\\
    &=(f(t),Au(t))_{L^2}-(B(u,u),Au)\nonumber\\
    &\leq \frac{1}{2\nu}\|f(t)\|^2+\frac{\nu}{2}\|Au(t)\|^2 + \|u\|_{L^{3}}\|\nabla u\|_{L^{2}}\|Au\|_{L^{6}} \nonumber\\
    &\leq \frac{1}{2\nu}\|f(t)\|^2+\frac{\nu}{2}\|Au(t)\|^2 +c^{2}\|u\|^{\frac{1}{2}}_{L^{2}}\|\nabla u\|^{\frac{3}{2}}_{L^{2}}\|\nabla Au\|_{L^{2}}\nonumber\\
    &\leq \frac{1}{2\nu}\|f(t)\|^2+\frac{\nu}{2}\|Au(t)\|^2 +\frac{c^{4}}{2\alpha^{2}\nu}\|u\|_{L^{2}}\|\nabla u\|^{3}_{L^{2}}+\frac{\alpha^{2}\nu}{2}\|\nabla Au\|^{2}_{L^{2}} 
    \nonumber\\&\leq \frac{1}{2\nu}\|f(t)\|^2+\frac{\nu}{2}\|Au(t)\|^2 +\frac{c^{4}M_{1}^{4}(t)}{2\alpha^{5}\nu}+\frac{\alpha^{2}\nu}{2}\|\nabla Au\|^{2}_{L^{2}}.
\end{align}

Therefore 
\begin{equation}\label{above4}
\frac{d}{dt}(\|\nabla u(t)\|^2+\alpha^2\|A u(t)\|^2)+\nu \lambda_{1}(\|\nabla u(t)\|^2+\alpha^2\| Au(t)\|^2)\leq
\frac{1}{\nu}\|f(t)\|^2+\frac{c^{4}M_{1}^{4}(t)}{\alpha^{5}\nu}  .  
\end{equation}

By  classical Gronwall's inequality applied in (\ref{above4}) (see \cite{evans2022partial}), we get
\begin{align}
    &\|\nabla u(t)\|^2+\alpha^2\|Au(t)\|^2\leq e^{-\nu\lambda_1t}(\|\nabla u(0)\|^2+\alpha^2\|A u(0)\|^2) \nonumber\\&+\int_0^te^{\nu\lambda_1(s-t)}\left(\frac{1}{\nu}\|f(s)\|^2+\frac{c^{4}M_{1}^{4}(s)}{\alpha^{5}\nu}\right)ds.\nonumber
    \end{align}
   We also have  
   \begin{align}
    &\int_0^te^{\nu\lambda_1(s-t)}\left(\frac{1}{\nu}\|f(s)\|^2+\frac{c^{4}M_{1}^{4}(s)}{\alpha^{5}\nu}\right)ds\leq \frac{1}{\nu^{2}\lambda_{1}}\sup_{s\geq 0}\|f(s)\|^2\nonumber\\
    &+\frac{2c^{4}}{\alpha^{5}\nu}\int_0^te^{\nu\lambda_1(s-t)}\left(e^{-2\nu\lambda_{1}s}\left(M^{2}_{A}+\alpha^{2}_{1}M^{2}_{B}\right)^{2}+\frac{1}{\lambda_1^4\nu^4}\cdot \sup\limits_{r\geq 0}\|f(r)\|^4\right)ds\nonumber\\
    &\leq
    \frac{1}{\nu^{2}\lambda_{1}}\sup_{s\geq 0}\|f(s)\|^2+\frac{2c^{4}}{\alpha_0^{5}\nu^{6}\lambda_1^5}\sup_{s\geq 0}\|f(s)\|^4\nonumber\\
    &+\frac{2c^{4}}{\alpha_0^{5}\nu^{2}\lambda_{1}}e^{-\nu\lambda_{1}t}\left(M^{2}_{A}+\alpha^{2}_{1}M^{2}_{B}\right)^{2}.\nonumber
\end{align}
Then we have (\ref{ine011}).

\end{proof}



\begin{lemma}\label{lem3}
    Let $u$ be a solution of \eqref{Bardina1Leray} and
let $M_{3}(t)$ be defined in (\ref{M3}), we have
    \begin{equation}
    \|u_t(s)\|\,\leq M_{3}(s),\,\, \forall\,\, s>0. \label{lemmm3}
    \end{equation}
\end{lemma}


\begin{proof}
     Applying the inverse operator $(I+\alpha^2A)^{-1}$ to \eqref{Bardina1Leray}, we obtain
$$u_t+\nu Au+(I+\alpha^2A)^{-1}B(u,u)=(I+\alpha^2A)^{-1}f.$$
Using (\ref{Poincare}), (\ref{Gagliardo-Nirenberg}) and (\ref{inequality-1}), we have
\begin{align}
   \|u_t\|\leq & \nu\|Au\|+\|(I+\alpha^2A)^{-1}B(u,u)\|+\|(I+\alpha^2A)^{-1}f\|\nonumber\\
    \leq &  \nu\|Au\|+\frac{1}{\alpha_{0}^2}\|B(u,u)\|_{D'}+\|f\|\nonumber\\
    \leq &  \nu\|Au\|+\frac{1}{\alpha^{2}_{0}\lambda_{1}^{\frac{1}{2}}}\|u^{2}\|+\|f\|\nonumber\\
    = &  \nu\|Au\|+\frac{1}{\alpha^{2}_{0}\lambda_{1}^{\frac{1}{2}}}\|u\|^{2}_{L^{4}}+\|f\|\nonumber\\
    \leq &  \nu\|Au\|+\frac{c^{2}}{\alpha^{2}_{0}\lambda_{1}^{\frac{1}{2}}}\|u\|^{\frac{1}{2}}\|\nabla u\|^{\frac{3}{2}}+\|f\|\nonumber\\
    \leq &  \nu\|Au\|+\frac{c^{2}}{\alpha^{2}_{0}\lambda_{1}^{\frac{3}{4}}}\|\nabla u\|^{2}+\|f\|.\nonumber
\end{align}
From Lemmas \ref{lem1} and \ref{lem12}, we obtain (\ref{lemmm3}).
\end{proof}



\begin{proposition}\label{teoBardina1}
 Let $u$ and $w_{n}$ be solutions of \eqref{Bardina1Leray} and \eqref{Bardina1LerayCDA}, respectively. Assume that the parameters $\eta_{n}$ and $N_{n}$ are chosen sufficiently large so that condition (\ref{estimates0002}) holds, and let $M^{2}_{(4,n)}(t,\eta_{n})$ be defined by (\ref{M4}). Then, the following estimate is valid  for all $t\in[t_{n},t_{n+1}]$:
 \begin{align}\label{error0}
  \displaystyle  \|w_{n}(t)\|^2 + \beta_{n}^2\|\nabla w_{n}(t)\|^2 \leq &  \displaystyle M^{2}_{(4,n)}(t,\eta_{n}).
\end{align}

\end{proposition}


\begin{proof}
Taking the $D^{'}$-dual action with $w_{n}$ in \eqref{Bardina1LerayCDA} and using \eqref{zero}, we have
\begin{align}
\displaystyle &\frac{1}{2}\frac{d}{dt}(\|w_{n}\|^{2}+\beta_{n}^2\|\nabla w_{n}\|^{2})+\nu \|\nabla w_{n}\|^{2}+\beta_{n}^2\nu\|Aw_{n}\|^{2}= (f,w_{n})-\eta_{n}(P_{N_{n}}w_{n}-w_{n},w_{n})\nonumber\\&-\eta_{n} \|w_{n}\|^{2}-\eta_{n}\beta_{n}^2 ( P_{N_{n}}w_{n}-w_{n},A w_{n})-\eta_{n}\beta_{n}^2 \|\nabla w_{n}\|^{2}
+\eta_n( (I+\beta_{n}^{2}A)P_{N_{n}}u, w_{n}).
\end{align}

We estimate each term of the right-hand side above using Young inequality and (\ref{Ih}):\\
$$ \|f\|\,\|w_{n}\|\leq \frac{1}{\eta_{n}}\|f\|^2+\frac{\eta_{n}}{4}\|w_{n}\|^2;$$
$$ \eta_{n}\|P_{N_{n}}w_{n}-w_{n}\|\,\|w_{n}\|\leq \eta_{n}\|P_{N_{n}}w_{n}-w_{n}\|^2+\frac{\eta_n}{4}\|w_{n}\|^2\leq\ds\frac{\eta_{n}\lambda_{1}^{-1}}{N^{2}_{n}}\|\nabla w_{n}\|^2+\frac{\eta_{n}}{4}\|w_{n}\|^2;$$
\begin{align} \eta_{n}\beta^{2}_{n}\|P_{N_{n}}w_{n}-w_{n}\|\,\|Aw_{n}\|&\leq \frac{\eta^{2}_{n}\beta^{2}_{n}}{2\nu}\|P_{N_{n}}w_{n}-w_{n}\|^2+\frac{\nu\beta^{2}_{n}}{2}\|Aw_{n}\|^2\nonumber\\&\nonumber\leq\ds\frac{\eta^{2}_{n}\beta^{2}_{n}\lambda_{1}^{-1}}{2\nu N^{2}_{n}}\|\nabla w_{n}\|^2+\frac{\nu\beta^{2}_{n}}{2}\|Aw_{n}\|^2;
\end{align}
$$ \eta_{n}\| P_{N_{n}}u\|\| w_{n}\| \leq \eta_{n}\| P_{N_{n}}u\|^{2}+\frac{\eta_{n}}{4}\|w_{n}\|^{2}; $$
$$ \eta_{n}\beta_{n}^{2}\|\nabla P_{N_{n}}u\|\|\nabla w_{n}\| \leq \frac{\eta_{n}\beta^{2}_{n}}{2}\|\nabla P_{N_{n}}u\|^{2}+\frac{\eta_{n}\beta^{2}_{n}}{2}\|\nabla w_{n}\|^{2}. $$

Therefore, we have
\begin{align}
\displaystyle \ds &\displaystyle\frac{d}{dt}(\|w_{n}\|^2  +\beta_{n}^2\|\nabla w_{n}\|^2)  +2\nu\|\nabla w_{n}\|^2+ \beta_{n}^2\nu\|Aw_{n}\|^2+ \nonumber\\
&\leq  \frac{2}{\eta_{n}}\|f\|^2+2\eta_{n}\| u\|^{2}+\eta_{n}\beta^{2}_{n}\|\nabla u\|^{2}-\frac{\eta_{n}}{2}\|w_{n}\|^{2}- \eta_{n}\beta^{2}_{n}\|\nabla w_{n}\|^{2}\nonumber\\&+\frac{\eta^{2}_{n}\beta^{2}_{n}\lambda_{1}^{-1}}{\nu N^{2}_{n}}\|\nabla w_{n}\|^2+\frac{2\eta_{n}\lambda_{1}^{-1}}{N^{2}_{n}}\|\nabla w_{n}\|^2.
\end{align}
By Lemma \ref{lem1}, we obtain 
\begin{align}
\displaystyle \ds \displaystyle\frac{d}{dt}(\|w_{n}\|^2  &+\beta_{n}^2\|\nabla w_{n}\|^2) +\frac{\eta_{n}}{2}\|w_{n}\|^{2} + \beta_{n}^2\nu\|Aw_{n}\|^2 \nonumber\\
 +&\left(2\nu +\eta_{n}\beta^{2}_{n}-\frac{\eta^{2}_{n}\beta^{2}_{n}\lambda_{1}^{-1}}{\nu N^{2}_{n}}-\frac{2\eta_{n}\lambda_{1}^{-1}}{N^{2}_{n}}\right)\|\nabla w_{n}\|^2\nonumber\\
\leq & \frac{2}{\eta_{n}}\sup_{s\geq 0}\|f(s)\|^2+\eta_{n}\left(2+\frac{\beta_{n}^{2}}{\alpha^{2}}\right)M^{2}_{1}(t)\nonumber.
\end{align}
From hypothesis (\ref{estimates0002}), we have $$2\nu +\eta_{n}\beta^{2}_{n}-\frac{\eta^{2}_{n}\beta^{2}_{n}\lambda_{1}^{-1}}{\nu N^{2}_{n}}-\frac{2\eta_{n}\lambda_{1}^{-1}}{N^{2}_{n}}\geq 2\nu +\eta_{n}\beta^{2}_{n}-\frac{\eta_{n}\beta^{2}_{n}}{2}-\nu>\frac{\eta_{n}\beta^{2}_{n}}{2}.$$
Then
\begin{align}
\displaystyle \ds \displaystyle\frac{d}{dt}(\|w_{n}\|^2  &+\beta_{n}^2\|\nabla w_{n}\|^2) +\frac{\eta_{n}}{2}\left(\|w_{n}\|^2  +\beta_{n}^2\|\nabla w_{n}\|^2\right) \nonumber\\
\leq & \frac{2}{\eta_{n}}\sup_{s\geq 0}\|f(s)\|^2+\eta_{n}\left(2+\frac{\beta_{n}^{2}}{\alpha^{2}}\right)M^{2}_{1}(t)\nonumber.
\end{align} 
Applying Gronwall’s inequality, we obtain
\begin{align}
&\|w_{n}(t)\|^2  +\beta_{n}^2\|\nabla w_{n}(t)\|^2  
\leq  e^{-\frac{\eta_{n}}{2}(t-t_{n})}\left(\|w_{n}(t_{n})\|^2  +\beta_{n}^2\|\nabla w_{n}(t_{n})\|^2\right)  \nonumber\\&+\int_{t_{n}}^{t}e^{\frac{\eta_{n}}{2}(s-t)}\left[\frac{2}{\eta_{n}}\sup_{r\geq 0}\|f(r)\|^2+\eta_{n}\left(2+\frac{\beta_{n}^{2}}{\alpha^{2}}\right)M^{2}_{1}(s)\right]ds
\nonumber\\
&\leq  e^{-\frac{\eta_{n}}{2}(t-t_{n})}\left(\|w_{n}(t_{n})\|^2  +\beta_{n}^2\|\nabla w_{n}(t_{n})\|^2\right)  \nonumber\\&+\frac{4}{\eta^{2}_{n}}\sup_{r\geq 0}\|f(r)\|^2+ 2\left(2+\frac{\beta_{n}^{2}}{\alpha^{2}}\right)M^{2}_{1}(t_{n}),
\nonumber
\end{align} 
and thus the estimate \eqref{error0} is obtained.

\end{proof}




\begin{proposition}\label{teoBardina2}
 Let $u$ and $w_{n}$ be solutions of \eqref{Bardina1Leray} and \eqref{Bardina1LerayCDA}, respectively. Assume that the parameters $\eta_{n}$ and $N_{n}$ are chosen sufficiently large so that conditions (\ref{estimates0001})-(\ref{estimates0002}) hold. Then, for all $t\in[t_{n},t_{n+1}]$, the following inequality holds for the difference between the physical and assimilated solutions, i.e., $g_{n}(t):=w_{n}(t)-u(t)$:
 \begin{align}\label{error}
  \displaystyle  \|g_{n}(t)\|^2 + \beta_{n}^2\|\nabla g_{n}(t)\|^2 \leq &   e^{\displaystyle  -\frac{\eta_{n}}{2}(t-t_{n})}\left(\|g_{n}(t_{n})\|^2  +\beta_{n}^2\|\nabla g_{n}(t_{n})\|^2\right)
 \nonumber\\&+ 4\left(\ds\frac{1}{\nu}M^{2}_{3}(t_{n})+\ds\frac{\nu}{\alpha_{0}^{2}}M^{2}_{2}(t_{n})\right)\frac{|\beta_{n}^2-\alpha^2|^2}{\eta_{n}\beta_{n}^2}.
\end{align}

\end{proposition}


\begin{proof}

Subtracting \eqref{Bardina1Leray} from \eqref{Bardina1LerayCDA} yields
\begin{align}\label{eqdifference11}
\displaystyle \frac{d}{dt}\big(g_{n}&+\beta_{n}^2Ag_{n}+(\beta_{n}^2-\alpha^2)Au\big)+\nu A(g_{n}+\beta_{n}^2Ag_{n}+(\beta_{n}^2-\alpha^2)Au)\nonumber \\
&+B(w_{n},w_{n})-B(u,u)=-\eta_{n} P_{N_{n}}g_{n}-\eta_{n}\beta_{n}^2 A P_{N_{n}}g_{n}, 
\end{align}
with $\nabla\cdot g_{n}=0$. Taking the $D^{'}$-dual action with $g_{n}$ in \eqref{eqdifference11}, using integration by parts and \eqref{zero}, we get
\begin{align}
\dfrac{1}{2}\displaystyle \frac{d}{dt}(\|g_{n}\|^{2}&+\beta_{n}^2\|\nabla g_{n}\|^{2})+(\beta_{n}^2-\alpha^2)(u_{t},Ag_{n})+\nu \|\nabla g_{n}\|^{2}+\beta_{n}^2\nu\|Ag_{n}\|^{2}\nonumber\\
&+\nu(\beta_{n}^2-\alpha^2)(Au,Ag_{n})+(B(g_{n},u),g_{n})=-\eta_{n}(P_{N_{n}}g_{n}-g_{n},g_{n})-\eta_{n} \|g_{n}\|^{2}\nonumber\\
&-\eta_{n}\beta_{n}^2 ( P_{N_{n}}g_{n}-g_{n},A g_{n})-\eta_{n}\beta_{n}^2 \|\nabla g_{n}\|^{2}.
\end{align}

Using general Hölder's inequality, we obtain
\begin{align}\label{xxx}
    \frac{1}{2}\frac{d}{dt}(\|g_{n}\|^2+\beta_{n}^2\|\nabla g_{n}\|^2)&+\nu(\|\nabla g_{n}\|^2+\beta_{n}^2\|Ag_{n}\|^2) \leq |\alpha^2-\beta_{n}^2|\|u_t\|\,\|Ag_{n}\| \nonumber\\
     &+\nu|\alpha^2-\beta_{n}^2|\|Au\|\,\|Ag_{n}\|+\|g_{n}\|^2_{L^4}\|\nabla u\|\nonumber\\&+\eta_{n}\|P_{N_{n}}g_{n}-g_{n}\|\,\|g_{n}\|-\eta_{n}\|g_{n}\|^2\nonumber\\
     &+\eta_{n}\beta_{n}^2\|P_{N_{n}}g_{n}-g_{n}\|\,\|Ag_{n}\|-\eta_{n}\beta_{n}^2\|\nabla g_{n}\|^2.
    \end{align}
We now estimate each term on the right-hand side using Young inequality together with (\ref{Gagliardo-Nirenberg}) and (\ref{Ih}), so that part of these contributions can be absorbed into the dissipation term:\\
$$|\alpha^2-\beta_{n}^2|\|u_t\|\,\|Ag_{n}\|\leq\ds\frac{\nu\beta_{n}^2}{4}\|Ag_{n}\|^2+\ds\frac{|\alpha^2-\beta_{n}^2|^2}{\nu\beta_{n}^2}\|u_t\|^2; $$
$$\nu|\alpha^2-\beta_{n}^2|\|Au\|\,\|Ag_{n}\|\leq \ds\frac{\nu\beta_{n}^2}{4}\|Ag_{n}\|^2+\ds\frac{\nu}{\beta_{n}^2}\|Au\|^2|\alpha^2-\beta_{n}^2|^2;$$
$$\|g_{n}\|^2_{L^4}\|\nabla u\|\leq c^{2}\|g_{n}\|^{1/2}\|\nabla g_{n}\|^{3/2}\|\nabla u\|\leq \ds\frac{\nu}{2}\|\nabla g_{n}\|^2+\frac{27c^4}{32 \nu^{3}}\|g_{n}\|^2\|\nabla u\|^4;$$
$$ \eta_{n}\|P_{N_{n}}g_{n}-g_{n}\|\,\|g_{n}\|\leq \frac{\eta_{n}}{2}\|P_{N_{n}}g_{n}-g_{n}\|^2+\frac{\eta_{n}}{2}\|g_{n}\|^2\leq\ds\frac{\eta_{n}\lambda^{-1}_{1}}{2N^{2}_{n}}\|\nabla g_{n}\|^2+\frac{\eta_{n}}{2}\|g_{n}\|^2;$$
\begin{align}
\eta_{n}\beta^{2}_{n}\|P_{N_{n}}g_{n}-g_{n}\|\,\|Ag_{n}\|&\leq \frac{\eta^{2}_{n}\beta^{2}_{n}}{\nu}\|P_{N_{n}}g_{n}-g_{n}\|^2+\frac{\nu\beta^{2}_{n}}{4}\|Ag_{n}\|^2\nonumber\\&\nonumber\leq\ds\frac{\eta^{2}_{n}\beta^{2}_{n}\lambda^{-1}_{1}}{\nu N^{2}_{n}}\|\nabla g_{n}\|^2+\frac{\nu\beta^{2}_{n}}{4}\|Ag_{n}\|^2.
\end{align}
With all estimates above in \eqref{xxx}, we get
\begin{eqnarray}\label{chk}
\displaystyle \ds &\displaystyle\frac{d}{dt}(\|g_{n}\|^2  +\beta_{n}^2\|\nabla g_{n}\|^2)  +\left(2\eta_{n}\beta^{2}_{n}+\nu-\frac{\eta_{n}\lambda^{-1}_{1}}{N^{2}_{n}}-\frac{2\eta^{2}_{n}\beta^{2}_{n}\lambda^{-1}_{1}}{\nu N^{2}_{n}}\right)\|\nabla g_{n}\|^2+ \frac{\beta_{n}^2\nu}{2}\|Ag_{n}\|^2 \nonumber\\
&\leq  \left(\ds\frac{27c^4}{16\nu^{3}}\|\nabla u\|^4-\eta_{n}\right)\|g_{n}\|^2+2\left(\ds\frac{1}{\nu\beta_{n}^2}\|u_t\|^2+\ds\frac{\nu}{\beta_{n}^2}\|Au\|^2\right)|\beta_{n}^2-\alpha^2|^2.\nonumber\\
\end{eqnarray}
By applying Lemmas \ref{lem1},  \ref{lem12} and \ref{lem3} to the estimate above, we remove the dependence on the norms of the physical solution $u(t)$, and obtain
\begin{align}\label{chk2}
\displaystyle  \displaystyle\frac{d}{dt}\left(\|g_{n}\|^2  +\beta_{n}^2\|\nabla g_{n}\|^2\right) &+\left(\eta_{n}-\ds\frac{27c^4M_1^4(t_{n})}{16\nu^{3}\alpha^4}\right)\|g_{n}\|^2+ \nonumber\\
& +\left(2\eta_{n}\beta^{2}_{n}+\nu-\frac{\eta_{n}\lambda^{-1}_{1}}{N^{2}_{n}}-\frac{2\eta^{2}_{n}\beta^{2}_{n}\lambda^{-1}_{1}}{\nu N^{2}_{n}}\right)\|\nabla g_{n}\|^2 + \frac{\beta_{n}^2\nu}{2}\|Ag_{n}\|^2  \nonumber\\
&  \leq 2\left(\ds\frac{1}{\nu\beta_{n}^2}M^{2}_{3}(t)+\ds\frac{\nu}{\beta_{n}^2\alpha^{2}}M^{2}_{2}(t)\right)|\beta_{n}^2-\alpha^2|^2.\nonumber\\
\end{align}
From (\ref{estimates0001})-(\ref{estimates0002}), we have

$$\eta_{n}-\ds\frac{27c^4M_1^4(t_{n})}{16\nu^{3}\alpha^4}\geq \eta_{n}-\ds\frac{27c^4M_1^4(t_{n})}{16\nu^{3}\alpha_{0}^4}\geq \frac{\eta_{n}}{2},$$
and $$2\eta_{n}\beta^{2}_{n}+\nu-\frac{\eta_{n}\lambda^{-1}_{1}}{N^{2}_{n}}-\frac{2\eta^{2}_{n}\beta^{2}_{n}\lambda^{-1}_{1}}{\nu N^{2}_{n}}\geq\frac{\eta_{n}\beta^{2}_{n}}{2},$$
which yields
\begin{align}\label{chk21}
\displaystyle  \displaystyle\frac{d}{dt}\left(\|g_{n}\|^2  +\beta_{n}^2\|\nabla g_{n}\|^2\right) &+\frac{\eta_{n}}{2}\left(\|g_{n}\|^2 +\beta^{2}_{n}\|\nabla g_{n}\|^2\right) \nonumber\\
\leq & 2\left(\ds\frac{1}{\nu\beta_{n}^2}M^{2}_{3}(t)+\ds\frac{\nu}{\beta_{n}^2\alpha^{2}}M^{2}_{2}(t)\right)|\beta_{n}^2-\alpha^2|^2.\nonumber\\
\end{align}
Finally, applying Gronwall’s inequality yields the desired error estimate \eqref{error}.

\end{proof}


\begin{proposition}\label{teoBardina25}
 Let $u$ and $w_{n}$ be solutions of \eqref{Bardina1Leray} and \eqref{Bardina1LerayCDA} on $[t_n,t_{n+1}]$, respectively. Assume there exist $\hat{t}_{n}$ with $t_{n}\leq \hat{t}_{n}<t_{n+1}$ and $\delta_{n}>0$, with $\delta_{n}$ defined in (\ref{deltaN}). Then, we have
\begin{align}\label{dfhbb}
 \displaystyle    
 |\alpha^{2} - \beta^{2}_{n+1}| \leq
 \displaystyle \frac{4c^{2}\lambda^{-\frac{1}{4}}_{1}}{\tilde{\delta}_{n}^{\frac{1}{2}}N_{n}^{\frac{1}{2}}}\sup_{\hat{t}_{n}\leq s\leq t_{n+1}}\left\{(\|\nabla w_{n}(s)\|+\|\nabla u(s)\|)\|\nabla g_{n}(s)\|\right\},
\end{align}
where 
\begin{equation}\label{dfbf}
\displaystyle \tilde{\delta}_{n}=\frac{\delta_{n}}{t_{n+1}-\hat{t}_{n}},
\end{equation}
with $c$ is given in (\ref{Gagliardo-Nirenberg}), and $\beta_{n+1}$ denoting the updated parameter obtained via the recovery algorithm (\ref{algBeta}).
\end{proposition}


\begin{proof}
Let $g_n=w_n-u$. Subtracting \eqref{Bardina1Leray} from \eqref{Bardina1LerayCDA} yields
\begin{align}\label{eqdifference2}
\displaystyle \frac{d}{dt}\big(g_{n}+\beta_{n}^2Ag_{n}&+(\beta_{n}^2-\alpha^2)Au\big)+\nu A(g_{n}+\beta_{n}^2Ag_{n}+(\beta_{n}^2-\alpha^2)Au)+\nonumber \\
+& B(w_{n},w_{n}) -B(u,u)=-\eta_{n} P_{N_{n}}g_{n}-\eta_{n}\beta_{n}^2 A P_{N_{n}}g_{n}, 
\end{align}
with $\nabla\cdot g_{n}=0$. Applying $P_{N_{n}}$ and taking the $D^{'}$-dual action with $\displaystyle P_{N_{n}}(u_{t}+\nu Au)$ in \eqref{eqdifference2}, and then integrating by parts, we obtain
\begin{align}
   &\left(\frac{d}{dt}P_{N_{n}}g_{n},P_{N_{n}}(u_{t}+\nu Au)\right)+\beta_{n}^2\left(\frac{d}{dt}\nabla P_{N_{n}}g_{n},\nabla P_{N_{n}}(u_{t}+\nu Au)\right)+\nonumber\\& +(\beta_{n}^2-\alpha^2)\|\nabla P_{N_{n}}(u_{t}+\nu Au)\|^{2}+\nu(\nabla P_{N_{n}}g_{n},\nabla P_{N_{n}}(u_{t}+\nu Au))+\nonumber\\&+\nu\beta_{n}^2(\nabla AP_{N_{n}}g_{n},\nabla P_{N_{n}}(u_{t}+\nu Au)) 
   +(P_{N_{n}}B(w_{n},w_{n}),P_{N_{n}}(u_{t}+\nu Au))+\nonumber\\
    &-(P_{N_{n}}B(u,u),P_{N_{n}}(u_{t}+\nu Au))\nonumber \\&=-\eta_{n}(P_{N_{n}}g_{n},P_{N_{n}}(u_{t}+\nu Au))-\beta_{n}^2\eta_{n}(\nabla P_{N_{n}}g_{n},\nabla P_{N_{n}}(u_{t}+\nu Au)).\label{eqdifference11A}
    \end{align}

Integrating over the time interval $[\hat{t}_{n},t_{n+1}]$, we obtain
\begin{align} \label{update2}
  (\alpha^2-\beta_{n}^2) \delta_{n}=&\int_{\hat{t}_{n}}^{t_{n+1}}\left(\frac{d}{dt}P_{N_{n}}g_{n},P_{N_{n}}(u_{t}+\nu Au)\right)+\beta_{n}^2 \left(\frac{d}{dt}\nabla P_{N_{n}}g_{n},\nabla P_{N_{n}}(u_{t}+\nu Au)\right)\nonumber\\
  +& \nu (\nabla P_{N_{n}}g_{n},\nabla P_{N_{n}}(u_{t}+\nu Au))+\nu\beta_{n}^2 (\nabla AP_{N_{n}}g_{n},\nabla P_{N_{n}}(u_{t}+\nu Au))
   \nonumber\\
    +& (P_{N_{n}}B(w_{n},w_{n}),P_{N_{n}}(u_{t}+\nu Au))- (P_{N_{n}}B(u,u),P_{N_{n}}(u_{t}+\nu Au))\nonumber \\
    +&\eta_{n} (P_{N_{n}}g_{n},P_{N_{n}}(u_{t}+\nu Au))+\beta_{n}^2\eta_{n} (\nabla P_{N_{n}}g_{n},\nabla P_{N_{n}}(u_{t}+\nu Au))\, ds.
    \end{align}
For estimating the difference $\alpha^2-\beta_{n+1}^2$, we use the equality (\ref{update2}), the update formula (\ref{algBeta}), and also (\ref{eqeq22}). Thus
\begin{align}\label{abs0101}
    &\alpha^{2} - \beta^{2}_{n+1}=(\alpha^{2} -\beta_{n}^{2})+ (\beta_{n}^{2}-\beta^{2}_{n+1})\nonumber\\&=  \displaystyle  \delta^{-1}_{n}\int_{\hat{t}_{n}}^{t_{n+1}}(P_{N_{n}}B(w_{n},w_{n}),P_{N_{n}}(u_{t}+\nu Au))-(P_{N_{n}}B(u,u),P_{N_{n}}(u_{t}+\nu Au))
\nonumber\\&+
(P_{N_{n}}B(P_{N_{n}}g_{n},P_{N_{n}}g_{n}),P_{N_{n}}(u_{t}+\nu Au))
    -(P_{N_{n}}B(w_{n},P_{N_{n}}g_{n}),P_{N_{n}}(u_{t}+\nu Au))\nonumber\\&
    -(P_{N_{n}}B(P_{N_{n}}g_{n},w_{n}),P_{N_{n}}(u_{t}+\nu Au))ds.
\end{align}
Moreover, since
\begin{align}
     \displaystyle  P_{N_{n}}B(w_{n},w_{n})-P_{N_{n}}B(u,u)=P_{N_{n}}B(g_{n},w_{n})+P_{N_{n}}B(w_{n},g_{n}) -P_{N_{n}}B(g_{n},g_{n}),\nonumber
\end{align}
we have
\begin{align}\label{ult1001.0}
   & \alpha^{2} - \beta^{2}_{n+1}=(\alpha^{2} -\beta_{n}^{2})+ (\beta_{n}^{2}-\beta^{2}_{n+1})\nonumber\\&=  \displaystyle  \delta^{-1}_{n}\int_{\hat{t}_{n}}^{t_{n+1}}(P_{N_{n}}B(g_{n},w_{n}),P_{N_{n}}(u_{t}+\nu Au))+(P_{N_{n}}B(w_{n},g_{n}),P_{N_{n}}(u_{t}+\nu Au))\nonumber\\&-(P_{N_{n}}B(g_{n},g_{n}),P_{N_{n}}(u_{t}+\nu Au))
+
(P_{N_{n}}B(P_{N_{n}}g_{n},P_{N_{n}}g_{n}),P_{N_{n}}(u_{t}+\nu Au))
    \nonumber\\&-(P_{N_{n}}B(w_{n},P_{N_{n}}g_{n}),P_{N_{n}}(u_{t}+\nu Au))
    -(P_{N_{n}}B(P_{N_{n}}g_{n},w_{n}),P_{N_{n}}(u_{t}+\nu Au))ds
    \nonumber\\&=
    \delta^{-1}_{n}\int_{\hat{t}_{n}}^{t_{n+1}}(P_{N_{n}}B(g_{n}-P_{N_{n}}g_{n},w_{n}),P_{N_{n}}(u_{t}+\nu Au))+\nonumber\\&+(P_{N_{n}}B(w_{n},g_{n}-P_{N_{n}}g_{n}),P_{N_{n}}(u_{t}+\nu Au))\nonumber\\&-(P_{N_{n}}B(g_{n},g_{n}),P_{N_{n}}(u_{t}+\nu Au))+
(P_{N_{n}}B(P_{N_{n}}g_{n},P_{N_{n}}g_{n}),P_{N_{n}}(u_{t}+\nu Au))ds
\nonumber\\&=\delta^{-1}_{n}\int_{\hat{t}_{n}}^{t_{n+1}}
    (P_{N_{n}}B(g_{n}-P_{N_{n}}g_{n},w_{n}),P_{N_{n}}(u_{t}+\nu Au))+\nonumber\\&+(P_{N_{n}}B(w_{n},g_{n}-P_{N_{n}}g_{n}),P_{N_{n}}(u_{t}+\nu Au))+\nonumber\\&+(P_{N_{n}}B(P_{N_{n}}g_{n}-g_{n},P_{N_{n}}g_{n}),P_{N_{n}}(u_{t}+\nu Au))
+\nonumber\\&+
    (P_{N_{n}}B(g_{n},P_{N_{n}}g_{n}-g_{n}),P_{N_{n}}(u_{t}+\nu Au))\,ds.
\end{align}

Therefore, we have
\begin{align}\label{ult1001.0aaa}
    \alpha^{2} - \beta^{2}_{n+1}
=\delta^{-1}_{n}\int_{\hat{t}_{n}}^{t_{n+1}}&
    (P_{N_{n}}B(g_{n}-P_{N_{n}}g_{n},w_{n}),P_{N_{n}}(u_{t}+\nu Au))+\nonumber\\+&(P_{N_{n}}B(w_{n},g_{n}-P_{N_{n}}g_{n}),P_{N_{n}}(u_{t}+\nu Au))+\nonumber\\+&(P_{N_{n}}B(P_{N_{n}}g_{n}-g_{n},P_{N_{n}}g_{n}),P_{N_{n}}(u_{t}+\nu Au))+
\nonumber\\+&
    (P_{N_{n}}B(g_{n},P_{N_{n}}g_{n}-g_{n}),P_{N_{n}}(u_{t}+\nu Au))\,ds.
\end{align}

We estimate each term of the right-hand side above using Young inequality, Gabliardo-Nirember inequalities (\ref{Gagliardo-Nirenberg}) and (\ref{Ih}) to get
\begin{align} \label{ult1001.1}
 \displaystyle    
 &|\alpha^{2} - \beta^{2}_{n+1}| \leq\displaystyle 2\delta^{-1}_{n}\int_{\hat{t}_{n}}^{t_{n+1}}(\|w_{n}\|_{L^{6}}+\|g_{n}\|_{L^{6}})\|g_{n}-P_{N_{n}}g_{n}\|_{L^{3}} \|\nabla P_{N_{n}}(u_{t}+\nu Au)\|\,ds \nonumber\\&\leq
 \displaystyle 4c^{2}\delta^{-1}_{n}\int_{\hat{t}_{n}}^{t_{n+1}}(\|\nabla w_{n}\|+\|\nabla u\|)\|g_{n}-P_{N_{n}}g_{n}\|^{\frac{1}{2}}\|\nabla g_{n}\|^{\frac{1}{2}} \|\nabla P_{N_{n}}(u_{t}+\nu Au)\|\,ds
 \nonumber\\&\leq
 \displaystyle \frac{4c^{2}\lambda^{-\frac{1}{4}}_{1}}{\delta_{n}N_{n}^{\frac{1}{2}}}\int_{\hat{t}_{n}}^{t_{n+1}}(\|\nabla w_{n}\|+\|\nabla u\|)\|\nabla g_{n}\|\|\nabla P_{N_{n}}(u_{t}+\nu Au)\|\,ds
 \nonumber\\&\leq
 \displaystyle \frac{4c^{2}\lambda^{-\frac{1}{4}}_{1}}{\delta_{n}N_{n}^{\frac{1}{2}}}\sup_{\hat{t}_{n}\leq s\leq t_{n+1}}\left\{(\|\nabla w_{n}(s)\|+\|\nabla u(s)\|)\|\nabla g_{n}(s)\|\right\}\int_{\hat{t}_{n}}^{t_{n+1}}\|\nabla P_{N_{n}}(u_{t}+\nu Au)\|\,ds
 \nonumber\\&\leq
 \displaystyle \frac{4c^{2}\lambda^{-\frac{1}{4}}_{1}}{\delta_{n}N_{n}^{\frac{1}{2}}}\sup_{\hat{t}_{n}\leq s\leq t_{n+1}}\left\{(\|\nabla w_{n}(s)\|+\|\nabla u(s)\|)\|\nabla g_{n}(s)\|\right\}\delta_{n}^{\frac{1}{2}}(t_{n+1}-\hat{t}_{n})^{\frac{1}{2}}
 \nonumber\\&\leq
 \displaystyle \frac{4c^{2}\lambda^{-\frac{1}{4}}_{1}}{\tilde{\delta}_{n}^{\frac{1}{2}}N_{n}^{\frac{1}{2}}}\sup_{\hat{t}_{n}\leq s\leq t_{n+1}}\left\{(\|\nabla w_{n}(s)\|+\|\nabla u(s)\|)\|\nabla g_{n}(s)\|\right\}.
\end{align}

Therefore, we conclude (\ref{dfhbb}).
\end{proof}

\section{Proof of the Main Result} 

\begin{proof}[Proof of Theorem \ref{theo121}]
\label{maintheoreticalresults3}
From the previous estimates (\ref{ine01}), (\ref{error0}), (\ref{error}) and (\ref{dfhbb}), we have
\begin{align}\label{dfbdfb}
 \displaystyle    
 &|\alpha^{2} - \beta^{2}_{n+1}| \leq
 \displaystyle \frac{4c^{2}\lambda^{-\frac{1}{4}}_{1}}{\tilde{\delta}_{n}^{\frac{1}{2}}N_{n}^{\frac{1}{2}}}\sup_{\hat{t}_{n}\leq s\leq t_{n+1}}\left\{(\|\nabla u(s)\|+\|\nabla w_{n}(s)\|)\|\nabla g_{n}(s)\|\right\}
 \nonumber\\
 &\leq
 \displaystyle \frac{4c^{2}\lambda^{-\frac{1}{4}}_{1}}{\tilde{\delta}_{n}^{\frac{1}{2}}N_{n}^{\frac{1}{2}}}\left[ \frac{M_{1}(\hat{t}_{n})}{\alpha_{0}}+\frac{M_{(4,n)}(\hat{t}_{n}, \eta_{n})}{\beta_{n}}\right]\cdot\frac{1}{\beta_n}
 \nonumber\\&
 \cdot\left[ e^{\displaystyle  -\frac{\eta_{n}}{4}(\hat{t}_{n}-t_{n})}\left(\|g_{n}(t_{n})\|  +\beta_{n}\|\nabla g_{n}(t_{n})\|\right)+ 2\left(\frac{\nu^{\frac{1}{2}}M_{2}(t_{n})}{\alpha_{0}}+\frac{M_{3}(t_{n})}{\nu^{\frac{1}{2}}}\right)\frac{|\beta_{n}^2-\alpha^2|}{\eta^{\frac{1}{2}}_{n}\beta_{n}}\right].
\end{align}

Using assumption (\ref{estimates000}), along with the definition of $\zeta_{n}$ given in (\ref{zeta}) and $\tilde{\delta_{n}}$ given in (\ref{dfbf}), we obtain 
$\zeta_{n} \leq \tilde{\delta_{n}}$.
Therefore,
\begin{align}
 \displaystyle    
 &|\alpha^{2} - \beta^{2}_{n+1}| \leq
 \displaystyle \frac{4c^{2}\lambda^{-\frac{1}{4}}_{1}}{\zeta_{n}^{\frac{1}{2}}N_{n}^{\frac{1}{2}}}\left[ \frac{M_{1}(\hat{t}_{n})}{\alpha_{0}}+\frac{M_{(4,n)}(\hat{t}_{n}, \eta_{n})}{\beta_{n}}\right]\cdot\frac{1}{\beta_n}\nonumber\\&
 \cdot\left[ e^{\displaystyle  -\frac{\eta_{n}}{4}(\hat{t}_{n}-t_{n})}\left(\|g_{n}(t_{n})\|  +\beta_{n}\|\nabla g_{n}(t_{n})\|\right)
 + 2\left(\frac{\nu^{\frac{1}{2}}M_{2}(t_{n})}{\alpha_{0}}+\frac{M_{3}(t_{n})}{\nu^{\frac{1}{2}}}\right)\frac{|\beta_{n}^2-\alpha^2|}{\eta^{\frac{1}{2}}_{n}\beta_{n}}\right]. \label{ulla1}
\end{align}

Now, let us use the above estimate to first obtain inequality (\ref{rrr1}). In fact,
from the results given in Lemma \ref{lem1} and Proposition \ref{teoBardina1} to estimate $g_{n}(t_{n})=w_{n}(t_{n})-u(t_{n})$, we obtain the following estimate
\begin{align}
 \displaystyle &\frac{4c^{2}\lambda^{-\frac{1}{4}}_{1}}{\beta_n\zeta_{n}^{\frac{1}{2}}N_{n}^{\frac{1}{2}}}\left[ \frac{M_{1}(\hat{t}_{n})}{\alpha_{0}}+\frac{M_{(4,n)}(\hat{t}_{n}, \eta_{n})}{\beta_{n}}\right]
  e^{\displaystyle  -\frac{\eta_{n}}{4}(\hat{t}_{n}-t_{n})}\cdot\left(\|g_{n}(t_{n})\|  +\beta_{n}\|\nabla g_{n}(t_{n})\|\right)\nonumber\\
  &\leq\frac{4c^{2}\lambda^{-\frac{1}{4}}_{1}}{\beta_n\zeta_{n}^{\frac{1}{2}}N_{n}^{\frac{1}{2}}}\left[ \frac{M_{1}(\hat{t}_{n})}{\alpha_{0}}+\frac{M_{(4,n)}(\hat{t}_{n}, \eta_{n})}{\beta_{n}}\right]
  e^{\displaystyle  -\frac{\eta_{n}}{4}(\hat{t}_{n}-t_{n})}\nonumber\\&\cdot\left(\|w_{n}(t_{n})\| + \|u(t_{n})\|  +\beta_{n}\|\nabla w_{n}(t_{n})\|+ \beta_{n}\|\nabla u(t_{n})\|\right)\nonumber\\
  &\leq 
  \displaystyle \frac{8c^{2}\lambda^{-\frac{1}{4}}_{1}}{\beta_n\zeta_{n}^{\frac{1}{2}}N_{n}^{\frac{1}{2}}}\left[ \frac{M_{1}(\hat{t}_{n})}{\alpha_{0}}+\frac{M_{(4,n)}(\hat{t}_{n}, \eta_{n})}{\beta_{n}}\right]
  e^{\displaystyle  -\frac{\eta_{n}}{4}(\hat{t}_{n}-t_{n})}\left[\left(1+\frac{\beta_{n}}{\alpha_{0}}\right)M_{1}(t_{n}) +M_{4,n}(t_{n},\eta_{n})\right].
  \label{ulla12}
\end{align}
 From  (\ref{ulla1}) and (\ref{ulla12}), along with assumptions (\ref{cond2})  and (\ref{cond1}), we obtain the following estimates for all $n\geq 1$:
$$ \displaystyle    
 |\alpha^{2} - \beta^{2}_{n+1}| \leq \frac{\varepsilon}{2}+ \frac{\varepsilon|\alpha^{2} - \beta^{2}_{n}|}{4|\alpha^{2}_{1}-\alpha^{2}_{0}|\chi_{1}(n)+4\varepsilon}.$$

By mathematical induction and the fact that $|\alpha^{2} - \beta^{2}_{1}|\leq |\alpha_{1}^{2} - \alpha_{0}^{2}|$, we obtain (\ref{rrr1}). \\

Now, let us use again the estimate (\ref{ulla1}) to obtain inequalities (\ref{rrr2}) and (\ref{rrr3}). In fact, using hypothesis (\ref{cond2}) and  (\ref{cond0}),  with (\ref{ulla1}), we also obtain
\begin{align}
 \displaystyle    
 |\alpha^{2} - \beta^{2}_{n+1}|
 &\leq
 \displaystyle 
 \frac{1}{4\nu\lambda_{1}^{\frac{3}{4}}}\left(\|g_{n}(t_{n})\|  +\beta_{n}\|\nabla g_{n}(t_{n})\|\right)
 + \frac{\varepsilon|\alpha^{2} - \beta^{2}_{n}|}{4|\alpha^{2}_{1}-\alpha^{2}_{0}|\chi_{1}(n)+4\varepsilon}
 \nonumber\\&\leq \displaystyle 
 \frac{1}{4\nu\lambda_{1}^{\frac{3}{4}}}\left(\|g_{n}(t_{n})\|  +\beta_{n}\|\nabla g_{n}(t_{n})\|\right)
 + \frac{|\alpha^{2} - \beta^{2}_{n}|}{4}
 .\label{last1}
\end{align}

On the other hand, combining inequalities (\ref{rrr1}) and (\ref{error}) with assumptions (\ref{cond4}) and (\ref{cond3}) and using the fact that $w_{n+1}(t_{n+1})= w_{n}(t_{n+1})$, we obtain
\begin{align}
  \displaystyle  &\|g_{n+1}(t_{n+1})\| + \beta_{n+1}\|\nabla g_{n+1}(t_{n+1})\|=
  \|g_{n}(t_{n+1})\| + \beta_{n+1}\|\nabla g_{n}(t_{n+1})\|
  \nonumber \\&\leq \max\left\{1,\frac{\beta_{n+1}}{\beta_{n}}\right\}(
  \|g_{n}(t_{n+1})\| + \beta_{n}\|\nabla g_{n}(t_{n+1})\|)
\nonumber \\&\leq 2\max\left\{1,\frac{\sqrt{\epsilon}+\alpha_{1}}{\beta_{n}}\right\}\nonumber\\&\cdot \left[e^{\displaystyle  -\frac{\eta_{n}}{4}(t_{n+1}-t_{n})}\left(\|g_{n}(t_{n})\|  +\beta_{n}\|\nabla g_{n}(t_{n})\|\right)
 + 2\left(\ds\frac{1}{\nu^{\frac{1}{2}}}M_{3}(t_{n})+\ds\frac{\nu^{\frac{1}{2}}}{\alpha_{0}}M_{2}(t_{n})\right)\frac{|\beta_{n}^2-\alpha^2|}{\eta_{n}^{\frac{1}{2}}\beta_{n}}\right]\nonumber\\&\leq
   \displaystyle  \frac{\left(\|g_{n}(t_{n})\|  +\beta_{n}\|\nabla g_{n}(t_{n})\|\right)}{4}
 + \nu\lambda_{1}^{\frac{3}{4}}\frac{|\beta_{n}^2-\alpha^2|}{4}.\label{last2}
\end{align}
Hence from (\ref{last1}) and (\ref{last2}), we have the system
\begin{eqnarray}\label{ullltttmmm}
\left\{
\begin{array}{l}
|\beta^{2}_{n+1}-\alpha^{2} |
 \leq 
 \displaystyle 
 \frac{\|g_{n}(t_{n})\|  +\beta_{n}\|\nabla g_{n}(t_{n})\|}{4 \nu\lambda_{1}^{\frac{3}{4}}}
 + \frac{|\beta_{n}^2-\alpha^2|}{4}, \\ \\
 \displaystyle\frac{\|g_{n+1}(t_{n+1})\| + \beta_{n+1}\|\nabla g_{n+1}(t_{n+1})\|}{\nu\lambda_{1}^{\frac{3}{4}}}\leq \displaystyle  \frac{\|g_{n}(t_{n})\|  +\beta_{n}\|\nabla g_{n}(t_{n})\|}{4 \nu\lambda_{1}^{\frac{3}{4}}}
 + \frac{|\beta_{n}^2-\alpha^2|}{4}.\\
\end{array}\right.
\end{eqnarray}
By applying mathematical induction, we then obtain (\ref{rrr2}) and (\ref{rrr3}). 
\end{proof}

\section{Conclusions}
In this work, we developed a novel parameter recovery algorithm for the three-dimensional viscous simplified Bardina turbulence model. We focused on recovering the length-scale parameter $\alpha$, which plays a fundamental role both as a filter length scale in turbulence modeling and as a regularization parameter that improves analytical tractability. Within a data assimilation framework, we introduced an algorithm for the approximate parameter $\beta$ that converges toward the true value $\alpha$ under suitable conditions.

All results presented here are theoretical. Extending the methodology to computational implementation remains a challenging task, which we are currently working on and will present in our next paper.

\textbf{Acknowledgment}
\vspace{3pt}

Jing Tian's work is partially supported by the NSF LEAPS-MPS Grant $\#2316894$.





\begin{thebibliography}{99}
	
	\bibitem{ALBANEZ2025109073}
D.~A.~F.~Albanez, M.~J.~Benvenutti, S.~Little, and J.~Tian,
\newblock Parameter analysis in continuous data assimilation for various turbulence models,
\newblock \emph{Commun. Nonlinear Sci. Numer. Simul.} \textbf{151} (2025), 109073.

\bibitem{albanez2018continuous}
D.~A.~Albanez and M.~J.~Benvenutti,
\newblock Continuous data assimilation algorithm for simplified Bardina model,
\newblock \emph{Evol. Equ. Control Theory} \textbf{7} (2018), no.~1, 33--52.

\bibitem{azouani2014continuous}
A.~Azouani, E.~Olson, and E.~S.~Titi,
\newblock Continuous data assimilation using general interpolant observables,
\newblock \emph{J. Nonlinear Sci.} \textbf{24} (2014), no.~2, 277--304.

\bibitem{bardina1980improved}
J.~Bardina, J.~H.~Ferziger, and W.~C.~Reynolds,
\newblock Improved subgrid-scale models for large-eddy simulation,
\newblock in \emph{13th Fluid and Plasmadynamics Conference} (1980), 1357.

\bibitem{animikh}
A.~Biswas and J.~Hudson,
\newblock Determining the viscosity of the Navier--Stokes equations from observations of finitely many modes,
\newblock \emph{Inverse Problems} \textbf{39} (2023), no.~12, 125012.

\bibitem{cao2006global}
Y.~Cao, E.~M.~Lunasin, and E.~S.~Titi,
\newblock Global well-posedness of the three-dimensional viscous and inviscid simplified Bardina turbulence models,
\newblock \emph{Comm. Math. Sci.} \textbf{4} (2006), no.~4, 823--848.

\bibitem{carlson2020parameter}
E.~Carlson, J.~Hudson, and A.~Larios,
\newblock Parameter recovery for the 2 dimensional Navier--Stokes equations via continuous data assimilation,
\newblock \emph{SIAM J. Sci. Comput.} \textbf{42} (2020), no.~1, A250--A270.

\bibitem{carlsonHudson2022}
E.~Carlson \emph{et al.},
\newblock Dynamically learning the parameters of a chaotic system using partial observations,
\newblock \emph{Discrete Contin. Dyn. Syst.} \textbf{42} (2022), no.~8, 3809--3839.

\bibitem{chen1999camassa}
S.~Chen \emph{et al.},
\newblock The Camassa--Holm equations and turbulence,
\newblock \emph{Physica D} \textbf{133} (1999), no.~1--4, 49--65.

\bibitem{cheskidov2005leray}
A.~Cheskidov, D.~D.~Holm, E.~Olson, and E.~S.~Titi,
\newblock On a Leray--$\alpha$ model of turbulence,
\newblock \emph{Proc. R. Soc. A} \textbf{461} (2005), no.~2055, 629--649.

\bibitem{constantin1988navier}
P.~Constantin and C.~Foias,
\newblock \emph{Navier-Stokes Equations},
\newblock University of Chicago Press, 1988.

\bibitem{evans2022partial}
L.~C.~Evans,
\newblock \emph{Partial Differential Equations},
\newblock Grad. Stud. Math., vol.~19, American Mathematical Society, 2022.

\bibitem{PhysRevFluids.9.054602}
A.~Farhat, A.~Larios, V.~R.~Martinez, and J.~P.~Whitehead,
\newblock Identifying the body force from partial observations of a two-dimensional incompressible velocity field,
\newblock \emph{Phys. Rev. Fluids} \textbf{9} (2024), no.~5, 054602.

\bibitem{Foias_Manley_Rosa_Temam_2001}
C.~Foias, O.~Manley, R.~Rosa, and R.~Temam,
\newblock \emph{Navier-Stokes Equations and Turbulence},
\newblock Encyclopedia of Mathematics and its Applications, vol.~83, Cambridge University Press, 2001.

\bibitem{foias2002three}
C.~Foias, D.~D.~Holm, and E.~S.~Titi,
\newblock The three dimensional viscous Camassa--Holm equations, and their relation to the Navier--Stokes equations and turbulence theory,
\newblock \emph{J. Dynam. Differential Equations} \textbf{14} (2002), no.~1, 1--35.

\bibitem{FoiasProdi}
C.~Foias and G.~Prodi,
\newblock Sur le comportement global des solutions non-stationnaires des équations de Navier--Stokes en dimension $2$,
\newblock \emph{Rend. Semin. Mat. Univ. Padova} \textbf{39} (1967), 1--34.

\bibitem{IlyinLunasinTiti2006}
A.~A.~Ilyin, E.~M.~Lunasin, and E.~S.~Titi,
\newblock A modified-Leray--$\alpha$ subgrid scale model of turbulence,
\newblock \emph{Nonlinearity} \textbf{19} (2006), no.~4, 879.

\bibitem{LaytonLewandowski2006}
W.~Layton and R.~Lewandowski,
\newblock On a well-posed turbulence model,
\newblock \emph{Discrete Contin. Dyn. Syst. Ser. B} \textbf{6} (2006), no.~1, 111.

\bibitem{Martinez2022}
V.~R.~Martinez,
\newblock Convergence analysis of a viscosity parameter recovery algorithm for the 2D Navier--Stokes equations,
\newblock \emph{Nonlinearity} \textbf{35} (2022), no.~5, 2241.


\bibitem{MartinezMurriWhitehead2025}
V.~R.~Martinez, J.~Murri, and J.~P.~Whitehead,
\newblock Relaxation-based schemes for on-the-fly parameter estimation in dissipative dynamical systems,
\newblock \emph{Inverse Problems} \textbf{41} (2025), no.~5, 055001.

\bibitem{Martinez2024}
V.~R.~Martinez,
\newblock On the reconstruction of unknown driving forces from low-mode observations in the 2D Navier--Stokes equations,
\newblock \emph{Proc. Roy. Soc. Edinburgh Sect. A} (2024), 1--24.

\bibitem{NeweyWhiteheadCarlson2025}
J.~Newey, J.~P.~Whitehead, and E.~Carlson,
\newblock Model discovery on the fly using continuous data assimilation,
\newblock \emph{J. Comput. Phys.} (2025), 114121.

\bibitem{PachevWhiteheadMcQuarrie2022}
B.~Pachev, J.~P.~Whitehead, and S.~A.~McQuarrie,
\newblock Concurrent multiparameter learning demonstrated on the Kuramoto--Sivashinsky equation,
\newblock \emph{SIAM J. Sci. Comput.} \textbf{44} (2022), no.~5, A2974--A2990.

\bibitem{RobinsonPierre2003}
J.~C.~Robinson and C.~Pierre,
\newblock \emph{Infinite-Dimensional Dynamical Systems: An Introduction to Dissipative Parabolic PDEs and the Theory of Global Attractors},
\newblock Cambridge Texts in Applied Mathematics, Cambridge University Press, 2003.

\bibitem{Temam2024}
R.~Temam,
\newblock \emph{Navier--Stokes Equations: Theory and Numerical Analysis},
\newblock vol.~343, American Mathematical Society, 2024.

\bibitem{Temam1995}
R.~Temam,
\newblock \emph{Navier--Stokes Equations and Nonlinear Functional Analysis},
\newblock Society for Industrial and Applied Mathematics, 1995.

\bibitem{WangJinHuang2023}
W.~Wang, C.~Jin, and Y.~Huang,
\newblock Recovering critical parameter for nonlinear Allen--Cahn equation by fully discrete continuous data assimilation algorithms,
\newblock \emph{Inverse Problems} \textbf{40} (2023), no.~1, 015008.

    \end{thebibliography}
\end{document}